\documentclass[12pt]{extarticle}
\usepackage{amsmath, amsthm, amssymb, color}
\usepackage[colorlinks=true,linkcolor=blue,urlcolor=blue]{hyperref}
\usepackage{graphicx}
\usepackage{caption}
\usepackage{mathtools}
\usepackage{enumerate}
\usepackage{verbatim}
\usepackage{tikz,tikz-cd,tikz-3dplot}
\usepackage{amssymb}
\usetikzlibrary{matrix}
\usetikzlibrary{arrows}
\usepackage{algorithm}
\usepackage[noend]{algpseudocode}
\usepackage{caption}
\usepackage[normalem]{ulem}
\usepackage{subcaption}
\usepackage{multicol}
\oddsidemargin \evensidemargin
\usepackage{makecell}
\usepackage{array}
\usepackage{enumitem}
\newtheorem{theorem}{Theorem}[section]

\theoremstyle{definition}

\newtheorem{problem}[theorem]{Problem}
\newtheorem{remark}[theorem]{Remark}

\newtheorem{examplex}[theorem]{Example}
\newenvironment{example}
{\pushQED{\qed}\examplex}
{\popQED\endexamplex}

\makeatletter
\def\Ddots{\mathinner{\mkern1mu\raise\p@
\vbox{\kern7\p@\hbox{.}}\mkern2mu
\raise4\p@\hbox{.}\mkern2mu\raise7\p@\hbox{.}\mkern1mu}}
\makeatother

\newcommand{\PP}{\mathbb{P}}
\newcommand{\RR}{\mathbb{R}}
\newcommand{\QQ}{\mathbb{Q}}
\newcommand{\CC}{\mathbb{C} }
\newcommand{\ZZ}{\mathbb{Z}}

\title{\bf Tropical Implicitization Revisited}

\author{Kemal Rose,
  Bernd Sturmfels and Simon Telen}
\date{}
\begin{document}
\maketitle

\begin{abstract}
\noindent
Tropical implicitization means computing the tropicalization of a unirational variety from its parametrization. In the case of a hypersurface, this amounts to finding the Newton polytope of the implicit equation, without computing its coefficients. We present a new implementation of this procedure in \texttt{Oscar.jl}. It solves challenging instances, and can be used for classical implicitization as well.  We also develop implicitization in higher codimension via Chow forms, and we pose several open questions.
\end{abstract}

\section{Introduction} \label{sec:1}
Let $X \subset \mathbb{C}^n$ be a $d$-dimensional affine variety defined as the closure of the image of a map
\begin{equation} \label{eq:parametrization}
f : \mathbb{C}^d \dashrightarrow \mathbb{C}^n, \quad t \longmapsto \bigl(f_1(t), f_2(t),\ldots, f_n(t) \bigr).
\end{equation}
Here $t = (t_1, \ldots, t_d)$, and $f_1,f_2,\ldots,f_n \in \mathbb{C}(t)$ are rational functions. The problem of \emph{implicitization} asks for the defining polynomial equations of $X$ in the coordinates $x = (x_1, \ldots, x_n)$ on $\mathbb{C}^n$. 
When $f_1, \ldots, f_n \in \mathbb{Q}(t)$ have rational coefficients, these equations can be computed via \emph{symbolic elimination} \cite[Chapter 3, \S 3]{cox1997ideals}. More precisely, one eliminates $t_1, \ldots, t_d$ from 
\begin{equation} \label{eq:graph}
x_1 - f_1(t) \,= \, x_2 - f_2(t) \, = \, \cdots \, = \, x_n - f_n(t) \, = \, 0 
\end{equation}
using Gr\"obner basis or resultant techniques. Unfortunately, these methods run out of steam 
for larger instances. This has motivated the question whether we can obtain interesting partial information in cases where computing the ideal of $X$ is out of reach.

Tropical geometry \cite{maclagan2021introduction} replaces an algebraic variety by a polyhedral complex which encodes many of its geometric properties. A commonly used slogan is that this complex serves as a \emph{combinatorial shadow} of the original variety. The \emph{tropicalization} ${\rm trop}(X)$ of $X$ is a pure $d$-dimensional polyhedral fan in $\mathbb{R}^n$, satisfying a balancing condition. The task of \emph{tropical implicitization} \cite{ST, STY} is to compute ${\rm trop}(X)$ from the data in \eqref{eq:parametrization}. This was the goal in the paper \cite{sturmfels2008tropical}, which includes demonstrations of an implementation called \texttt{TrIm}. Theorem 6.4 in \cite{sturmfels2008tropical} suggests the following two-step procedure for performing tropical implicitization:
\begin{enumerate}
\item Compute the tropicalization of the graph $\Gamma_f$ of $f$, given by the $n$ equations in \eqref{eq:graph}.
The result is the $d$-dimensional balanced fan
${\rm trop}(\Gamma_f)$ in the product space $\,\mathbb{R}^d \times \mathbb{R}^n$.
\item Project this fan to $\mathbb{R}^n$ and assign appropriate multiplicities to each cone in the image. 
\end{enumerate}
This is illustrated in Figure \ref{fig: classical vs tropical implicitization}. Computing ${\rm trop}(\Gamma_f)$ in step 1 can be complicated in general. It involves the computation of a \emph{tropical basis} for the ideal generated by $x_1 - f_1(t), \ldots, x_n-f_n(t)$. 

In this article, we consider two different assumptions on the map $f$. Both assumptions circumvent the tropical basis computation, and are relevant in practice. First, in Section \ref{sec:2}, we assume that the functions $f_i$ are Laurent polynomials which are generic with respect to their Newton polytopes. This is the assumption in \cite{ST, STY, sturmfels2008tropical}. It reduces step 1 above to computing the stable intersection of $n$ 
codimension one fans in $\mathbb{R}^d \times \mathbb{R}^n$. Second, in Section~\ref{sec:3}, we assume that $f$ 
is the composition of
 a linear map $\lambda: \mathbb{C}^d \rightarrow \mathbb{C}^\ell$
followed by a Laurent monomial map  $\mu: \mathbb{C}^\ell \dashrightarrow\mathbb{C}^n$.
In symbols, we have  $f = \mu \circ \lambda$.  
 This allows to compute ${\rm trop}(X)$ as the linear projection of a tropical linear space in $\mathbb{R}^\ell$. 
 For details see \cite[Section 5.5]{maclagan2021introduction}.
 An important special case arises from the computation of \emph{tropical $A$-discriminants} \cite{dickenstein2007tropical}. 

Tropical implicitization is a first step towards classical implicitization. Let $X = V(F)$ be the hypersurface defined by a polynomial $F \in \mathbb{C}[x]$. Then ${\rm trop}(X)$ is the union of the $(n-1)$-dimensional cones in the normal fan of the Newton polytope ${\cal N}(F)$, decorated with multiplicities. From ${\rm trop}(X)$ we can recover ${\cal N}(F)$. The key ingredient  is a \emph{vertex oracle} which, for a generic weight vector $w \in \mathbb{R}^n$, returns the vertex $v$ of ${\cal N}(F)$ which minimizes the dot product with $w$ on ${\cal N}(F)$. The algorithm realizing the oracle is suggested by \cite[Theorem 2.2]{dickenstein2007tropical}.  We provide an implementation using \texttt{Oscar.jl}
and use it to recover ${\cal N}(F)$ via the algorithm in \cite{huggins}. Once we have the Newton polytope ${\cal N}(F)$, we can find $F$ via (numerical) linear algebra. The task is to compute the unique kernel vector of a matrix constructed via \emph{numerical integration} \cite{corless2001numerical} or \emph{sampling} \cite{breiding2018learning,emiris2013implicitization}. Sampling is preferred when the $f_i$ have rational coefficients. We can then use the parametrization to find rational points on $X$, and $F$ can be computed using exact arithmetic over $\mathbb{Q}$. However, the size of the matrix is the number of lattice points in ${\cal N}(F)$, and we may have to resort to floating point arithmetic when this number is too large. An alternative is rational reconstruction from linear algebra over finite fields. We discuss these techniques in Section \ref{sec:4}. We use them to solve instances for which elimination via Gr\"obner bases does not terminate within reasonable time. 

If the $f_i$ are Laurent polynomials which are generic with respect to their Newton polytopes, as in Section \ref{sec:2}, then ${\cal N}(F)$ is a \emph{mixed fiber polytope} \cite{emiris2007computing,esterov2008elimination,STY}. Our implementation in \texttt{Oscar.jl} for computing ${\cal N}(F)$ gives a practical way of computing mixed fiber polytopes.

When $\dim X < n-1$, we present a new way of finding its implicit equations
 from ${\rm trop}(X)$. This is the topic of Section \ref{sec:5}. The idea is to pass through the \emph{Chow form} ${\rm Ch}(X)$ of $X$ \cite{DS}. The polytope we compute is the \emph{Chow polytope} ${\cal C}(X)$, which is a linear projection of the Newton polytope ${\cal N}({\rm Ch}(X))$. This computation rests on a result by Fink \cite{Fink}, which describes the (weighted) normal fan of ${\cal C}(X)$ in terms of ${\rm trop}(X)$. We explain how to recover ${\rm Ch}(X)$ from ${\cal C}(X)$, using the parameterizing functions and an appropriate ansatz. Defining equations for $X$ are obtained from ${\rm Ch}(X)$ in the standard manner \cite[Proposition 3.1]{DS}. 

The implementation of the algorithms supporting this work have benefited from the flexibility provided by \texttt{Oscar.jl}. The possibility to combine polyhedral  computations with symbolic linear and nonlinear algebra in the same environment has greatly simplified the task. This feature has been our incentive to revisit tropical implicitization. Throughout the article, we include several open problems and computational challenges which we hope will inspire the reader to join this effort.
Our software and data 
are made available in the {\tt MathRepo} collection at MPI-MiS via
 \url{https://mathrepo.mis.mpg.de/TropicalImplicitization}.

\section{Generic tropical implicitization} \label{sec:2}

In this section, we start with $n$ Laurent polynomials in $n$ variables with complex coefficients:
\[ \qquad f_i \, =\, \sum_{a \in A_i} c_{i,a} \, t^a \,\in \,
\mathbb{C}[t_1^{\pm 1}, \ldots, t_d^{\pm 1}]
\quad {\rm for} \quad i \, = \, 1,2,\ldots, n. \]
We use these Laurent polynomials in \eqref{eq:parametrization}. 
The tuple $f = (f_1, \ldots, f_n)$ gives a map $(\mathbb{C}^*)^d \rightarrow (\mathbb{C}^*)^n$.
Let $X \subset (\mathbb{C}^*)^n$ be
 the closure of the image of $f$.  Our first task is to find its tropicalization ${\rm trop}(X)$. 
In Section \ref{sec:4}, we use ${\rm trop}(X)$ for classical implicitization.  As a~set, 
\[ {\rm trop}(X) \, = \, \{ w \in \mathbb{R}^n \,:\, {\rm in}_w(I(X)) \text{ does not contain a monomial } \}\, \, \subset \mathbb{R}^n. \]
Here $I(X) \subset \mathbb{C}[x_1^{\pm1}, \ldots, x_n^{\pm 1}]$ is the vanishing ideal of $X$, and ${\rm in}_w$ takes the initial ideal with respect to the weight vector $w$. It is well known that ${\rm trop}(X)$ is the support of a fan $\Sigma$ of dimension ${\rm dim}(X)$. This fan is not unique, but for the purposes of this text we can choose any fan $\Sigma$ with support ${\rm trop}(X)$. Assigning a \emph{multiplicity} $m_\sigma$ to each top dimensional cone $\sigma \in \Sigma$ in the appropriate way \cite[Definition 3.4.3]{maclagan2021introduction}, the fan $\Sigma$ is \emph{balanced} \cite[Theorem 3.4.14]{maclagan2021introduction}. We will see that these multiplicities are crucial when using ${\rm trop}(X)$ for implicitization. 

Classically, the variety $X \subset (\mathbb{C}^*)^n$ is the closure of the projection  $\Gamma_f \rightarrow (\mathbb{C}^*)^n$ of the~graph
\[ \Gamma_f \, =\,\{(x,t) \in (\mathbb{C}^*)^n \times (\mathbb{C}^*)^{d} \, : \, x_1 - f_1(t) = 0, \ldots, x_n - f_n(t) = 0 \} \]
onto the $n$ $x$-coordinates. It turns out this has an easy tropical analog. 
\begin{theorem} \label{thm:project}
Let $X = \overline{{\rm im} \, f} \subset (\mathbb{C}^*)^n$. The tropical variety ${\rm trop}(X)$ is the image of the projection ${\rm trop}(\Gamma_f) \rightarrow \mathbb{R}^n$, where ${\rm trop}(\Gamma_f) \subset \mathbb{R}^n \times \mathbb{R}^d$ is the tropicalization of the graph of $f$. 
\end{theorem}
This is an instance of \cite[Theorem 2.1]{STY}. See also \cite[Theorem 6.4]{sturmfels2008tropical}. We can thus obtain ${\rm trop}(X)$ from ${\rm trop}(\Gamma_f)$ via a simple projection. 
However, Theorem \ref{thm:project} is only useful in practice when ${\rm trop}(\Gamma_f)$ is easy to compute. Our next theorem describes ${\rm trop}(\Gamma_f)$ under the assumption that the $f_i$ are \emph{generic} with respect to their Newton polytopes ${\cal N}(f_i)$. It uses the following notation.
For a polytope $P \subset \mathbb{R}^k$ and a vector~$w \in (\mathbb{R}^k)^*$, we write $P_w = \{ p \in P \, : \, w \cdot p \leq w \cdot q \text{ for all } q \in P \}$. In words, $P_w$ is the face of $P$ supported by $w$.

 \begin{theorem} \label{thm:tropgraph}
Suppose $f_i$  is generic with respect to ${\cal N}(f_i)$, and let $P_i = {\cal N}(x_i - f_i(t)) \subset \mathbb{R}^n \times \mathbb{R}^d$ for $i = 1, \ldots, n$. 
The tropical variety $\,{\rm trop}(\Gamma_f)$ is the support of a $d$-dimensional subfan of the normal fan of $P= P_1 + \cdots + P_n$. It consists of the normal cones $\sigma$ of $P$ for which the face polytopes $(P_1)_w,\ldots,(P_n)_w$ have positive mixed volume ${\rm MV}_w$ in the affine lattice of $P_w$, for each $w \in {\rm int}(\sigma)$. Moreover, the multiplicity $m_\sigma$ of $\sigma$ in ${\rm trop}(\Gamma_f)$ equals~${\rm MV}_w$.
\end{theorem}
This is \cite[Theorem 4.3]{ST}. We illustrate this theorem for a parametric plane curve.

\begin{figure}[h]
\centering
\includegraphics[width = 15cm]{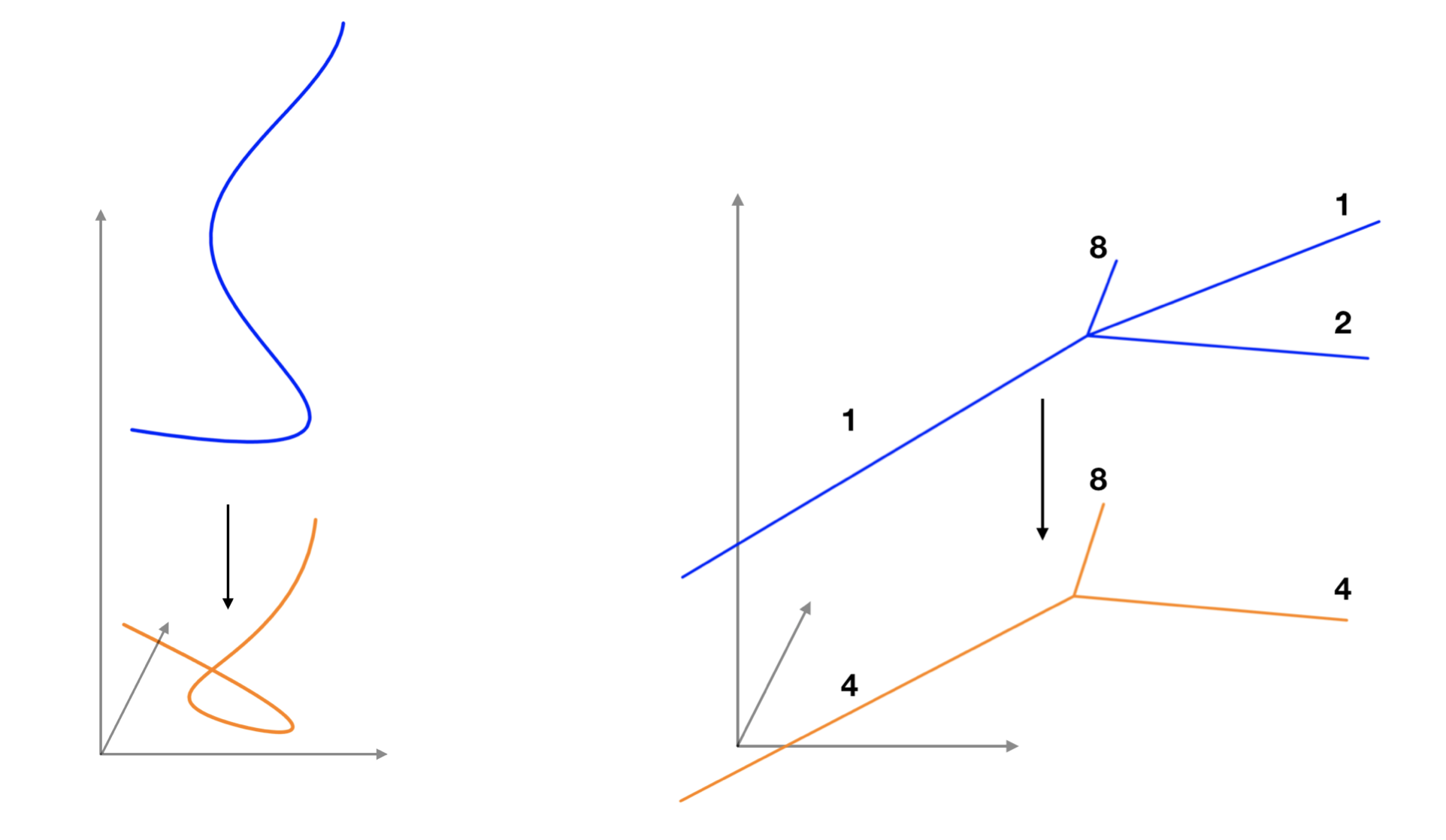} \vspace{-0.14in}
\caption{Classical (left) and tropical (right) implicitization of a parametric plane curve. }
\label{fig: classical vs tropical implicitization}
\end{figure}

\begin{example}
\label{ex: tropical curve}
Consider the parametrization $f = (f_1, f_2):  \CC^* \longrightarrow ( \CC^*)^2$ given by
\[
f_1 = 11 \, t^2 + 5 \, t^3 - t^4   \quad \text{and} \quad 
f_2 =  11 + 11 \, t + 7 \, t^8.
\]
The image is the plane curve $C = \overline{\operatorname{im}f}$ given by the implicit equation $F(x,y) = 0$, with 
\begin{equation} \label{eq:Fcurve}
F \, = \,  2401 \, x^8 
- 1372 \,x^6y
- 422576\,x^5y +
 \cdots + y^4 +  \cdots + 1247565503668. 
\end{equation}
This has $25$ terms, one for each lattice point of ${\cal N}(F)$, shown on the right side of Figure \ref{fig:polytopes}. 
The Newton polytopes of $x-f_1 $ and $y - f_2 $ of $\Gamma_f$ are the triangles seen on the left side.

\begin{figure}[h]
\centering
\includegraphics[width = 13cm]{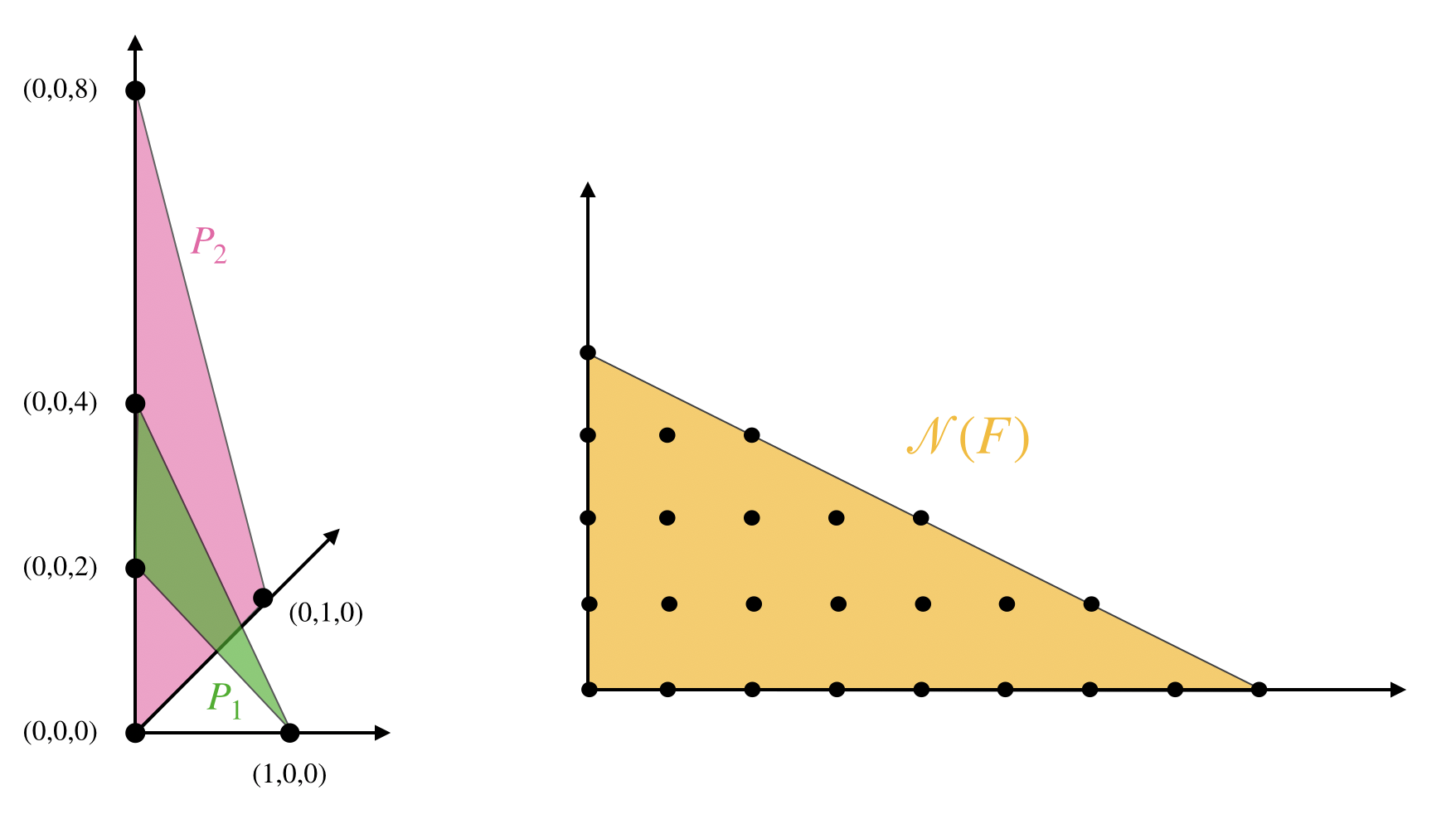} \vspace{-0.12in}
\caption{Newton polytopes $P_1, P_2$ and ${\cal N}(F)$ from Example \ref{ex: tropical curve}.} 
\label{fig:polytopes}
\end{figure}

The tropical curve ${\rm trop}(\Gamma_f)$
can be constructed according to Theorem \ref{thm:tropgraph}.
It is shown in blue on the right of Figure \ref{fig: classical vs tropical implicitization}. The result is a balanced, one-dimensional fan
with four rays:
\[
{\rm trop}(\Gamma_f) \, \, = \, \,
\RR_+  \cdot (1, 0,0)  \,\, \cup \,  \,
\RR_+   \cdot (-4, -8,-1) \,\,  \cup \,  \,
\RR_+   \cdot (0, 1,0) \,\, \cup \,  \,
\RR_+  \cdot (2, 0,1) ,
\]
with respective multiplicities $2, 1, 8$ and $1$. 
We demonstrate how to obtain these multiplicities.

Consider the primitive ray generator $w =  (2, 0,1) $,  revealing the face polytopes
\[
 \left(P_1 \right)_w \,=\, \operatorname{conv}(    (0,0, 2), (1, 0,0)  )
 \quad \text{  and  }  \quad \left(P_2 \right)_w 
 \,=\,  \operatorname{conv}(   (0, 0, 0), (0, 1, 0)   ).
\]
The multiplicity of  $\RR_+  \cdot (2,0,1)$ in ${\rm trop}(\Gamma_f)$ can be computed as the  mixed volume of the line segments $\left(P_1 \right)_w$ and $\left(P_2 \right)_w$ inside the
$2$-dimensional lattice define by the affine hull of their Minkowski sum.
This is the mixed volume ${\rm MV}_w$, and we find that it is equal to $1$.
\end{example}

Now that we know ${\rm trop}(\Gamma_f)$ and its multiplicities (when the Laurent polynomials $f_i$ are generic), and we know that ${\rm trop}(X)$ is obtained from its projection, it remains to determine the multiplicities of ${\rm trop}(X)$ from those of ${\rm trop}(\Gamma_f)$. The answer is given by \cite[Theorem 1.1]{ST}, which is the second part of \cite[Theorem 6.4]{sturmfels2008tropical}. In order to recall the formula, we introduce some more notation. Let $\Sigma_X$ be a fan in $\mathbb{R}^n$ whose support is ${\rm trop}(X)$, and $\Sigma_{\Gamma_f}$ a fan in $\mathbb{R}^n \times \mathbb{R}^d$ whose support is ${\rm trop}(\Gamma_f)$. Let $v$ be a point in the interior of a top dimensional cone $\sigma_v \in \Sigma_X$. We write $\mathbb{L}_v$ for the linear span of a small open neighborhood of $v$ in ${\rm trop}(X)$. Similarly, $w \in {\rm int}(\sigma_w)$ for a top dimensional cone $\sigma_w \in \Sigma_{\Gamma_f}$ defines a linear space $\mathbb{L}_w$. If the projection $\Gamma_f \rightarrow X$ is generically finite of degree $\delta$, then the multiplicity of $\sigma_v \in \Sigma_X$ is 
\begin{equation}
\label{eq: pushforward}
 m_{\sigma_v} \, = \, \frac{1}{\delta} \sum_{w \in \pi^{-1}(v)} m_{\sigma_w} \cdot \text{index}
 \bigl(\,\mathbb{L}_v \cap \mathbb{Z}^n : \pi(\mathbb{L}_w \cap \mathbb{Z}^{n + d})\, \bigr). 
\end{equation}
Here $\pi$ is the projection $\mathbb{R}^n \times \mathbb{R}^d \rightarrow \mathbb{R}^n$ and the sum is over all points $w$ in the pre-image of $v$ under 
the map $\pi_{|{\rm trop}(\Gamma_f)} : {\rm trop}(\Gamma_f) \rightarrow {\rm trop}(X)$. It is assumed that there are only finitely many such points, and each of them lies in the interior of a top dimensional cone of $\Sigma_{\Gamma_f}$. 

With the choice of weights (\ref{eq: pushforward}), the image fan is balanced. 
This is a non-trivial fact,  derived in a more general setting in
 \cite[Lemma 3.6.3]{maclagan2021introduction}.
 See also \cite[Theorem 6.5.16]{maclagan2021introduction}
 for a textbook discussion of tropical implicitization in the context of
 geometric tropicalization.

\begin{example} \label{ex:tropcurve}
According to Theorem  \ref{thm:project},
the  tropical curve ${\rm trop}(\Gamma_f)$  projects to
${\rm trop}(C) $.
This is displayed on the right of 
Figure \ref{fig: classical vs tropical implicitization},
where ${\rm trop}(C) $ is shown in orange as the fan 
\[ {\rm trop}(C) \,\, = \, \, \RR_+  \cdot  (1, 0) \,\,  \cup \, \,
\RR_+  \cdot (-1,-2) \, \, \cup \,  \, \RR_+   \cdot (0, 1). \]
This fan is balanced with ray multiplicities $ 4, 4, 8$, in that order.
We demonstrate the computation of $m_\rho = 4$ for the first ray $\rho =  \RR_+  \cdot  (1, 0)$ using \eqref{eq: pushforward}. The ray $\hat{\rho} = \RR_+  \cdot  (2, 0,1)$ of ${\rm Trop}(\Gamma_f)$ projects to $\rho$. Its primitive ray generator $(2,0,1)$ projects to the imprimitive lattice vector $(2, 0)$.
The contribution of $\hat{\rho} = \RR_+ \cdot (2, 0,1)$ to the multiplicity $m_\rho$ is the product of two numbers: its intrinsic multiplicity $m_{\hat{\rho}} = 1$, and the lattice index $2$.
The ray $\RR_+ \cdot( 1, 0, 0)$ also projects to $\rho$, which leads to a total of $m_\rho \, = \, 1\cdot 2 + 2 \cdot 1$. 
The tropical curve
${\rm trop}(C) $ equals the normal fan of the Newton polytope
$\mathcal{N}(F)$, displayed on the right of Figure \ref{fig:polytopes}.
\end{example}

The discussion above leads to Algorithm \ref{alg:gentropimp}, which makes the results in this section effective.
It
 takes the Newton polytopes $Q_i = {\cal N}(f_i)$ as an input, and returns
the tropicalization of $X = \overline{\operatorname{im}f}$. Here the Laurent polynomials $f_i$ are assumed to be generic with respect to their Newton polytopes $Q_i$. 
The output is  a set of pairs
$(m_\tau, \tau)$, where $\tau \subset \RR^n$ is a cone, and $m_\tau$ is a positive integer. The tropical hypersurface
${\rm trop}(X)$ is the union of all these cones $\tau$, and the multiplicity of ${\rm trop}(X)$ at a generic point $x$ is the sum
$\sum_{x\in \tau} m_\tau$.

We warn the reader that, although the union of all cones
$\tau$ forms the support of a fan, the collection of cones itself is generally not a fan.
This representation of a tropical variety is unconventional. However, it is easy to compute and convenient for our algorithmic purposes. 

\begin{algorithm} 
\caption{Generic tropical implicitization}
\begin{algorithmic}[1] 
\Procedure{getTropicalCycle($Q_1, \ldots, Q_n$)}{}
\For{ $i  \in  \{1, \dots, n\} $}
	    \State $P_i \gets {\rm conv}(e_{i} \cup Q_i)$
	    \label{line: def P_i}
\EndFor
\State $P \gets P_1 + \cdots +P_n$
\State $\Sigma \gets \text{normal fan of }P$
\State ${\rm trop}(X) \gets  \emptyset $
\For{ $\sigma \in \Sigma$}
	    \State $ m_\sigma \gets {\rm MV}( (P_1)_\sigma, \ldots, (P_n)_\sigma) $
	    \label{line: cpt MV}
\If{$m_\sigma > 0$}
\State $\tau \gets \pi(\sigma)$
\State $m_{\rm lattice}  \gets {\rm index}(\mathbb{L}_\tau \cap \mathbb{Z}^n: \pi(\mathbb{L}_\sigma \cap \mathbb{Z}^{n+d}))$
\label{line: lattice mult}
\State ${\rm trop}(X) \gets {\rm trop}(X)  \cup \{( m_{\rm lattice} \cdot m_\sigma, \tau)  \}  $
\EndIf
\EndFor
\State \Return \textit{${\rm trop}(X)$}
\EndProcedure
\end{algorithmic}
\label{alg:gentropimp}
\end{algorithm}

We now explain
 Algorithm \ref{alg:gentropimp}.
The polytopes $P_1, \dots, P_n$ in line \ref{line: def P_i} are the Newton polytopes of
the equations $x_1 -f_1, \dots, x_n-f_n$ of $\Gamma_f$. 
The standard basis $e_1,e_2,\ldots\,$ of $\mathbb{R}^{n +d}$ is indexed by the variables $x_1, \ldots, x_n, t_1, \ldots, t_d$ in that order. 
Following Theorem \ref{thm:tropgraph}, 
Algorithm \ref{alg:gentropimp} selects all cones $\sigma$ in the normal fan of $P_1+ \cdots + P_n$ that contribute to the tropicalization 
${\rm trop}(\Gamma_f)$.  Line \ref{line: cpt MV} computes the mixed volume
$m_\sigma = {\rm MV}( (P_1)_\sigma, \ldots, (P_n)_\sigma) $, where $(P_i)_\sigma = (P_i)_w$ for any $w \in {\rm int}(\sigma)$.
We denote by $\pi(\sigma)$ the projection of $\sigma \subset \RR^{n+d}$ to the first $n$ coordinates.
The multiplicity with which $\pi(\sigma)$ contributes to ${\rm trop}(X)$
is computed in lines  \ref{line: cpt MV} and \ref{line: lattice mult}. Based on \eqref{eq: pushforward},
it is the product of $m_\sigma$ with  
 the index of the lattice
 $\pi(\mathbb{L}_\sigma \cap \mathbb{Z}^{n+d})$ in the lattice 
$\mathbb{L}_\tau \cap \mathbb{Z}^n$.
Here $\mathbb{L}_\sigma$
is the linear span of $\sigma$ and
$\mathbb{L}_\tau = \pi(\mathbb{L}_\sigma)$. We implemented Algorithm \ref{alg:gentropimp} in {\tt Julia}.

\begin{example}
\label{example: code: cycle frome parametrization}
We show how to apply our {\tt Julia}  implementation to Example \ref{ex: tropical curve}:
\begin{verbatim}
using TropicalImplicitization, Oscar
R, (t,) = polynomial_ring(QQ,["t"])
f1 = 11*t^2 + 5*t^3 - 1*t^4 
f2 = 11 + 11*t + 7*t^8
Q1 = newton_polytope(f1)
Q2 = newton_polytope(f2)
newton_pols = [Q1, Q2]
cone_list, weight_list = get_tropical_cycle(newton_pols)
\end{verbatim}
The lists \texttt{cone\_list} and \texttt{weight\_list} returned by our program  
have four elements each. The first list contains the planar cones
\[
\RR_+  \cdot  (1, 0) \, , \quad
\RR_+  \cdot (1,0) \, ,  \quad  
\RR_+  \cdot (0,1) \, , \quad
\RR_+  \cdot (-1,-2),
\]
and the second list consists of their respective multiplicities $(2,2,8,4)$. Notice that $\RR_+  \cdot  (1, 0)$ appears twice, and its multiplicity is split up as $4 = 2+2$, like in Example \ref{ex:tropcurve}.
\end{example}

\begin{problem}
Suppose the coefficients of $f_1,\ldots,f_n$ lie in a 
field with a non-trivial valuation, such as the
$p$-adic numbers $\mathbb{Q}_p$ or the
Puiseux series $\CC \{ \! \{ t \} \! \}$. While the theory
of tropical implicitization generalizes nicely to this
setting, with balanced fans replaced by balanced polyhedral complexes,
useful algorithms and their implementations are yet to be developed.
 \end{problem}

\section{$A$-discriminants} \label{sec:3}

We fix a $d \times n$ integer matrix $A$ of rank $d$ which has the
vector  $(1,1,\ldots,1)$ in its row span. 
The associated $(d-1)$-dimensional projective toric variety $X_A$
is the closure in $\PP^{n-1}$ of the~set
\begin{equation}
  \label{eq:XA}
  \bigl\{ \,(t^{a_1} : t^{a_2} : \cdots : t^{a_n}) \in \PP^{n-1} \,:\,
  t = (t_1,\ldots,t_d) \in (\CC^*)^d \,\bigr\} .
\end{equation}
Here $a_i$ denotes the $i$th column of the matrix $A$.
We are interested in the dual variety $X_A^*$, which parametrizes
hyperplanes that are tangent to $X_A$ at some points. Equivalently,
$X_A^*$ is the closure in $\PP^{n-1}$ of the set of
points $x = (x_1:x_2 : \cdots : x_n)$ such that the
hypersurface
\begin{equation}
\label{eq:hypersurface} \bigl \{ \, t \in (\CC^*)^d \,: \, \sum_{i=1}^n x_i \,t^{a_i} \,=\, 0 \,\bigr\}
\end{equation}
has a singular point. The variety $X_A^*$ is irreducible, and it is usually 
a hypersurface. The {\em $A$-discriminant}
 $\Delta_A$ is the
unique (up to scaling) irreducible polynomial vanishing on $X_A^*$.

In this section we address the following computational problem:
given the matrix $A$, compute its $A$-discriminant $\Delta_A$.
Along the way, we will discover whether $X_A^*$ is not a hypersurface.
In this event, we turn to Section 5, and we compute its Chow form instead.

Our algorithm is based on the {\em Horn uniformization}, which writes
$X_A^*$ as the image of a map whose coordinates are products of linear 
forms. We follow the exposition given in \cite{dickenstein2007tropical}.
For additional information, see the book references
in \cite[Section 9.3.F]{GKZ} and \cite[Section 5.5]{maclagan2021introduction}.
Given two vectors $u$ and $v$ in $(\CC^*)^n$, we define
$\, u \star v  = (u_1 v_1: u_2 v_2 : \cdots : u_n v_n) \,\in \, \PP^{n-1}$.
If $U$ and $V$ are varieties in $\PP^{n-1}$, neither contained in a
coordinate hyperplane, then their {\em Hadamard product}
$U \star V$ is the closure of all such points $ u \star v$, 
where $u \in U$ and $v \in V$.

\begin{theorem}[Horn Uniformization]
\label{thm: Horn Uniformization}
The dual variety $X_A^*$ is the Hadamard product in
$\PP^{n-1}$ of the
$(d-1)$-dimensional toric variety $X_A$ with
an $(n-d-1)$-dimensional linear space:
\begin{equation}
\label{eq:hornformula} X_A^* \,\,= \,\,X_A \,\star \, {\rm kernel}(A) . 
\end{equation}
\end{theorem}

We illustrate this theorem with several examples.
In each of them, we refer to 
the $(d-1)$-dimensional polytope $Q = {\rm conv}(a_1,a_2,\ldots,a_n)$,
and we fix an  $(n-d) \times n$-matrix $B$
whose rows span the kernel of $A$. In polytope language,
$B$ is a Gale transform of the polytope $Q$.
For (\ref{eq:hornformula}),
we introduce unknowns $u = (u_1,\ldots,u_{n-d})$ and we write
$u B$ for vectors in ${\rm kernel}(A)$.

\begin{example}[Determinant]
Fix $n = k^2$ and $d=2k-1$, for some integer $k \geq 2$,
and let $A$ represent the linear map that extracts
the row sums and column sums of a $k \times k$ matrix. Naively, this matrix has $2k$ rows, but only $2k-1$ of them are linearly independent. Here
$Q = \Delta_{k-1} \times \Delta_{k-1}$ is the product of two
$(k-1)$-simplices. The toric variety $X_A$ consists of
$k \times k$ matrices of rank $1$ and
$X_A^*$ consists of $k \times k$ matrices of rank $\leq k-1$.
We parametrize $X_A^*$  by  the Hadamard product of a
rank $1$ matrix with a matrix whose row and columns are zero.
E.g.,~for $k=3$, 
the Horn uniformization writes all singular $3 \times 3$ matrices as follows:
\begin{equation}
\label{eq:horn33}
\begin{bmatrix} 
t_1 t_4 \, u_1 & t_1 t_5 (u_2-u_1)  &  t_1 t_6 (-u_2)  \\
t_2 t_4 (u_3-u_1) & t_2 t_5 (u_1 -u_2-u_3+u_4)   &  t_2 t_6 (u_2-u_4) \\
t_3 t_4 (-u_3)              & t_3 t_5 (u_3-u_4)                             &  t_3 t_6 \,u_4
\end{bmatrix}.
\end{equation}
This matrix has $(t_4^{-1},t_5^{-1}, t_6^{-1} )^t$ in its right kernel and 
$(t_1^{-1},t_2^{-1},t_3^{-1})$ in its left kernel.
The $A$-discriminant $\Delta_A$ is the determinant of a square matrix,
which obviously vanishes on (\ref{eq:horn33}).
\end{example}

\begin{example}[Resultant]
The resultant of a square system of homogeneous polynomials
is the $A$-discriminant where $A$ is the Cayley configuration 
of the given monomial supports. 
We examine the Sylvester resultant $\Delta_A$ of two binary quadrics ($d=3,n=6$). We set
$$ A \,\,= \,\, \begin{bmatrix}
\,1 \, & \, 1 \, & \, 1 \, & \, 0 \, & \, 0 \, & \, 0 \,\,\\
\,0 \, & \, 0 \, & \, 0 \, & \, 1 \, & \, 1 \, & \, 1 \,\,\\
\,0 \, & \, 1 \, & \, 2 \, & \, 0 \, & \, 1 \, & \, 2 \,\,\\
\end{bmatrix} \quad {\rm and} \quad
B \,\, = \,\, 
\begin{bmatrix}
\,1 & -2 & \phantom{-}1 & \phantom{-}0 & \phantom{-}0 & \phantom{-}0 \,\,\\
\,0 & \phantom{-}0 & \phantom{-}0 & \phantom{-} 1 & -2 & \phantom{-}1\,\,  \\
 \, 1 & -1 & \phantom{-}0 &  -1 & \phantom{-} 1 & \phantom{-}0\,\, \\
 \end{bmatrix}.
$$
This yields the following parametrization for pairs of univariate quadrics with a common~zero:
$$
\begin{matrix}
 x_1 \,+\, x_2 \,z + x_3 \,z^2  &=& t_1 (t_3 u_1 \,z - u_1 - u_3) (t_3 \,z-1),    \\
 x_4\, +\, x_5 \,z + x_6 \,z^2 &=& t_2 (t_3 u_2 \,z - u_2 - u_3) (t_3 \,z-1).
\end{matrix}
$$
These Horn uniformizations exist for resultants of polynomials in any number of variables.
\end{example}

\begin{example}[Hyperdeterminant]
\label{ex: hyperdeterminant}
The hyperdeterminant of a multidimensional tensor vanishes whenever the
hypersurface defined by the associated multilinear form is singular.
In our notation, this is the $A$-discriminant $\Delta_A$ where the columns of
$A$ are the vertices of a product of simplices. As an illustration, we here present the
Horn uniformization for the hyperdeterminant of format $2 \times 2 \times 2$.
Here $n = 8$ and our configuration is the regular $3$-cube:
$$ 
A \, = \,\begin{bmatrix}
\,1 & 1 & 1 & 1 & 1 & 1 & 1 & 1 \, \\
\,0 & 0 & 0 & 0 & 1 & 1 & 1 & 1 \,\\
\,0 & 0 & 1 & 1 & 0 & 0 & 1 & 1 \,\\
\,0 & 1 & 0 & 1 & 0 & 1 & 0 & 1 \, \\
\end{bmatrix}
\quad {\rm and} \quad
B = \begin{bmatrix}
\,1 &  -1 &  -1 &  \phantom{-}1 & \phantom{-} 0 &  \phantom{-}0 & \phantom{-} 0 &  0\, \,\\
\,0 &  \phantom{-}0 &  \phantom{-}0 &  \phantom{-}0 & \phantom{-} 1 &  -1 &  -1 &  1 \,\,\\
\,1 &  -1 &   \phantom{-}0 &   \phantom{-}0 &  -1 &   \phantom{-}1 &   \phantom{-}0 &  0 \,\,\\
\,1 &   \phantom{-}0 &  -1 &   \phantom{-}0 &  -1 &   \phantom{-}0 &   \phantom{-}1 &  0 \,\,\\
\end{bmatrix}.
$$
These two matrices yield the following map
from $\CC^8$ into the space of $2 \times 2 \times 2$ tensors
$$ \begin{matrix}
x_{000}  =   t_1(u_1+u_3+u_4), & 
x_{001} =    t_1 t_4 (-u_1-u_3), &
x_{010} =    t_1 t_3 (-u_1-u_4), &
x_{011} =    t_1 t_3 t_4 u_1, \\
x_{100} =    t_1 t_2 (u_2-u_3-u_4) ,&
x_{101} =   t_1 t_2 t_4 (u_3-u_2), &
x_{110} =   t_1 t_2 t_3 (u_4-u_2), &
x_{111} =    t_1 t_2 t_3 t_4 u_2.
\end{matrix}
$$
Implicitization of this parametrization gives us the hyperdeterminant:
$$ \begin{matrix}
 \Delta_A \,\,=\,\,  x_{000}^2x_{111}^2
+x_{001}^2 x_{110}^2
+x_{011}^2 x_{100}^2
+x_{010}^2 x_{101}^2
+ 4 x_{000} x_{011} x_{101} x_{110}
+4 x_{001} x_{010} x_{100} x_{111} &
\\
-\,2 x_{000} x_{001} x_{110} x_{111}
-2 x_{000} x_{010} x_{101} x_{111}
-2 x_{000} x_{011} x_{100} x_{111}  & 
\\
-\,2 x_{001} x_{010} x_{101} x_{110}
-2 x_{001} x_{011} x_{100}    x_{110}
-2 x_{010} x_{011} x_{100} x_{101}. & \qedhere
 \end{matrix}     
$$
\end{example}
\smallskip

We now return to tropical implicitization. Our aim is to compute the tropical variety
${\rm trop}(X_A^*)$ directly from $A$.  Here we identify $X_A^*$ with its affine cone in
$(\CC^*)^n$. If $X_A^*$ has codimension $1$ then ${\rm trop}(X_A^*)$
 is an $(n-1)$-dimensional balanced fan in $\mathbb{R}^n$, with a one-dimensional lineality space.
This is the normal fan of the Newton polytope of the $A$-discriminant $\Delta_A$.
We recover the polytope from the fan using 
Algorithm \ref{alg: vertex_oracle} below; see also
\cite[Remark 3.3.11]{maclagan2021introduction}.

The Horn uniformization of Theorem \ref{thm: Horn Uniformization} gives a convenient way of computing ${\rm trop}(X_A^*)$. 
It is an instance of parametrizations given by
monomials in linear forms. These
admit an elegant solution to the tropical implicitization problem; see
\cite[Section 5.5]{maclagan2021introduction}.
Let $U$ and $V$ be integer matrices of size $r \times m$ and $s \times r$ respectively. The rows of $V$ are $v_1, \dots, v_s \in \mathbb{Z}^r$.
We denote by $\lambda_U$ the linear map defined by $U$, and by $\mu_V$
the monomial map specified by~$V$:
\begin{multicols}{2}
\noindent
\begin{align*}
\lambda_U : (\CC^*)^m &\dashrightarrow (\CC^*)^r \\
v &\longmapsto U\, v
\end{align*}
\begin{align*}
\mu_V : (\CC^*)^r &\longrightarrow (\CC^*)^s \\
x &\longmapsto ( x^{v_1}, \dots, x^{v_s} ).
\end{align*}
\end{multicols}

The composition of these maps gives the unirational variety
$Y_{U,V} = \overline{ \text{im}( \mu_V \circ \lambda_U )}$ in $(\CC^*)^s$.
Its tropicalization ${\rm trop}(Y_{U,V})$ is obtained by tropicalizing the map $\mu_V \circ \lambda_U$.
We begin with the tropical linear space ${\rm trop}({\rm im}\, \lambda_U)$. This is computed
purely combinatorially, as  the \emph{Bergman fan} of the matroid of $U$;
see \cite[Section 4.2]{maclagan2021introduction}.
 The monomial map $\mu_V$ tropicalizes to the linear map $V: \mathbb{R}^r \rightarrow \mathbb{R}^s$. The following result is  \cite[Theorem 3.1]{dickenstein2007tropical} and
\cite[Theorem 5.5.1]{maclagan2021introduction}.

\begin{theorem}
The tropical variety ${\rm trop}(Y_{U,V})$
 is the image, as a balanced fan via \cite[Lemma 3.6.3]{maclagan2021introduction}, of the Bergman fan 
${\rm trop}({\rm im} \,\lambda_U)$ under the linear map $\RR^r \rightarrow \RR^s$ given by 
 $V$.
\end{theorem}

\noindent By Theorem \ref{thm: Horn Uniformization}, the affine cone over the $A$-discriminant in $(\mathbb{C}^*)^n$ is the variety
$Y_{U,V}$ with
\begin{equation}
\label{eq: Matrices U, V}
U = \begin{pmatrix}
B^t & 0 \\
0 & I_d \\
\end{pmatrix} , \quad
V = \begin{pmatrix}
I_n & A^t \\
\end{pmatrix} .
\end{equation}
Here $m=s=n$ and $r= n+d$.
This leads to Algorithm \ref{alg:tropAdisc} for computing
${\rm trop}(X_A^*)$.

\begin{algorithm}
\caption{Compute the tropical $A$-discriminant}
\begin{algorithmic}[1]
\Procedure{${\rm getTropADisc}$}{$A$}
\State $B \gets \text{Gale dual of }A$
\label{line: GaleDual}
\State $U \gets   \begin{pmatrix}
B^t & 0 \\
0 & I_d \\
\end{pmatrix} $
\State $V \gets \begin{pmatrix}
I_n & A^t \\
\end{pmatrix} $
\State $M \gets \text {matroid of } U$
\label{line: 1 tropical linear space}
\State ${\rm trop}({\rm im}\, \lambda_U) \gets  \text{Bergman fan of } M $
\label{line: 2 tropical linear space}
\State ${\rm trop}(X_A^*) \gets  \emptyset $
\For{$(m_\sigma, \sigma) \in {\rm trop}({\rm im}\, \lambda_U)$}
\label{line: project trop U}
	\State $\tau \gets V   \sigma$
	\State $m_{\rm lattice}  \gets {\rm index}(\mathbb{L}_\tau \cap \mathbb{Z}^n: V(\mathbb{L}_\sigma \cap \mathbb{Z}^{n+d}))$
	\State ${\rm trop}(X_A^*) \gets {\rm trop}(X_A^*) \cup \{(m_\sigma \cdot m_{\rm lattice}, \tau)\}$
\EndFor
\State \Return \textit{${\rm trop}(X_A^*) $}
\EndProcedure
\end{algorithmic}
\label{alg:tropAdisc}
\end{algorithm}

The matrix $B$ in line \ref{line: GaleDual} is  Gale dual to $A$.
Using the symbolic linear algebra functionality provided by \texttt{Oscar.jl}, 
we find this with the command
\texttt{nullspace(A)}.
Lines \ref{line: 1 tropical linear space} and~\ref{line: 2 tropical linear space}
compute the tropicalization ${\rm trop}({\rm im}\, \lambda_U)$ of the column span of $U$.
They are based on the \texttt{Oscar.jl} commands
\texttt{Oscar.Polymake.matroid.Matroid(VECTORS = U)}
and \texttt{Oscar.Polymake.tropical.matroid\_fan\{min\}(matroid)}.
From line \ref{line: project trop U} on, the algorithm 
computes a  projection of the Bergman fan ${\rm trop}({\rm im}\, \lambda_U)$.
This is analogous to Algorithm \ref{alg:gentropimp}.

\begin{example}
\label{example: code:  hyperdeterminant}
We compute the tropicalized $2 \times 2 \times 2$ hyperdeterminant from Example 
\ref{ex: hyperdeterminant}:
\begin{verbatim}
A = [1 1 1 1 1 1 1 1; 0 0 0 0 1 1 1 1; 0 0 1 1 0 0 1 1; 0 1 0 1 0 1 0 1]
cone_list, weight_list = get_trop_A_disc(A)
\end{verbatim}
The result consists of 32 {7-dimensional} cones and a list of their multiplicities, constituting the weighted normal fan of the Newton polytope $\mathcal{N}(\Delta_A)$.
The following code 
uses an implementation of Algorithm~\ref{alg: vertex_oracle} below. It
computes $\mathcal{N}(\Delta_A)$, its lattice points, and its f-vector.
\begin{verbatim}
Delta = get_polytope_from_cycle(cone_list, weight_list)
f_vec, lattice_pts = f_vector(Delta), lattice_points(Delta)
\end{verbatim}
The result is $\mathtt{f\_vec} = (6,14,16,8)$, and \texttt{lattice\_pts} contains the 12 exponents of $\Delta_A$.
\end{example}

Mixed discriminants \cite{sandra} are special cases of $A$-discriminants.
We discuss a non-trivial~one.

\begin{example} \label{ex:43400}
We revisit \cite[Example 5.1]{dickenstein2007tropical}.
Here, $d=4$ and $n=8$, and we fix the matrix
$$ A \,= \, \begin{small} \begin{bmatrix}
1 & 1 & 1 & 1 & 0 & 0 & 0 & 0 \\
0 & 0 & 0 & 0 & 1 & 1 & 1 & 1 \\
2 & 3 & 5 & 7 & 11 & 13 & 17 & 19 \\
19 & 17 & 13 & 11 & 7 & 5 & 3 & 2 
\end{bmatrix}. \end{small}
$$
This represents the following sparse system of two polynomial equations in two variables:
$$ x_1 s^2 t^{19} + x_2 s^3 t^{17} + x_3 s^5 t^{13} + x_4 s^7 t^{11} 
\,\, = \,\, x_5 s^{11} t^7 + x_6 s^{13} t^5 + x_7 s^{17} t^3 + x_8 s^{19} t^2 \,\,=\,\, 0 .
$$
The mixed volume of the two Newton polygons equals $39$, so we
expect $39$ common solutions in $(\CC^*)^2$.
The $A$-discriminant $\Delta_A$ is the condition for two of these solutions to come together.
We know from \cite[Example 5.1]{dickenstein2007tropical} that $\Delta_A$ is
a polynomial of degree $126$ in $x_1,\ldots,x_8$. 
Our software computes that
$\mathcal{N}(\Delta_A)$ has $43400$ lattice points, and f-vector $(45,92,63,16)$.
\end{example}
\begin{problem}
Computing the coefficients of $\Delta_A$ in Example \ref{ex:43400} amounts to solving a linear system of equations over $\mathbb{Q}$ with $43400$ unknowns. This is one of our topics in Section~\ref{sec:4}. Solving that system is hard for at least two reasons. First, systems of this size are beyond the reach of symbolic black box solvers on most personal computers at present. Second, the large condition number and the unbalanced nature of the coefficients of the implicit equation, as in \eqref{eq:Fcurve}, hinder the naive use of numerical linear algebra. It is an interesting challenge to develop symbolic or mixed symbolic-numerical techniques for solving such problems. 
\end{problem}

\section{Polytope reconstruction and interpolation} \label{sec:4}

Suppose that $X$ is an irreducible hypersurface in $(\CC^*)^n$, given by 
its parametrization (\ref{eq:parametrization}) or (\ref{eq:hornformula}).
Using Algorithms \ref{alg:gentropimp} and \ref{alg:tropAdisc},
we have computed the tropical variety $\Sigma_X = {\rm trop}(X)$.
Thus, $\Sigma_X$ is a balanced fan of dimension $d = n-1$ in $\RR^n$,
represented by a collection of weighted cones.
Our aim in this section is to compute the polynomial $F \in \mathbb{C}[x_1, \ldots, x_n]$ 
that defines the hypersurface $X$. We identify $X$ with its closure in $\CC^n$.
 This makes $F$  uniquely defined up to scaling. In particular, 
 the Newton polytope $P = \mathcal{N}(F)$ is uniquely specified.
  
Before spelling out the details, we summarize our approach. First, we compute the Newton polytope 
$P$ from the balanced fan $\Sigma_X$. This relies on Theorem \ref{thm:troptonewt} below. Second, 
we find the polynomial $F$ from the parametrization by interpolation. Here we use the ansatz 
\begin{equation} \label{eq:ansatzF}
 F(x) \quad = \, \sum_{a \, \in \, {\cal N}(F) \cap \mathbb{Z}^n} c_a \, x^a,
 \end{equation}
and we determine the unknown coefficients $c_a$ by evaluating \eqref{eq:ansatzF} at many points $x$ on $X$. 

We start with computing $P = {\cal N}(F)$. The fan $\Sigma_X$ is dual to the Newton polytope $P$,
namely,
 it is the $(n-1)$-skeleton of the normal fan of $P$. Taking into account the multiplicities of all 
 maximal cones of $\Sigma_X$, we can go back and forth between $P$ and $\Sigma_X$. 
Obtaining $\Sigma_X$ and its multiplicities from $P$ is straightforward. The multiplicity of an $(n-1)$-dimensional cone in $\Sigma_X$ is the lattice length of the corresponding edge of $P$. The other direction is more interesting to us: we want to compute $P$ from the output of Algorithm \ref{alg:gentropimp}.
This is discussed in \cite[Remark 3.3.11]{maclagan2021introduction}.
The main tool is a \emph{vertex oracle}, provided by \cite[Theorem~2.2]{dickenstein2007tropical}. 

\begin{theorem} \label{thm:troptonewt}
Let $X = V(F) \subset \mathbb{C}^n$ be a hypersurface, whose tropicalization ${\rm trop}(X)
\subset \RR^n$
  is the support of a fan $\Sigma_X$. For a generic weight vector $w \in \mathbb{R}^n$, the vertex ${\cal N}(F)_w$~is 
\[ \sum_{i=1}^n \left (\sum_{\sigma \in \Sigma_X} m_\sigma \cdot {\rm IM}(w + \mathbb{R}_+ \cdot e_i, \sigma) \right ) \cdot e_i.\]
 Here $e_i$ is a standard basis vector, and
 the inner sum is over all maximal cones of $\,\Sigma_X$.
 \end{theorem}
 
\begin{figure}[h]
\centering
\includegraphics[width = 15cm]{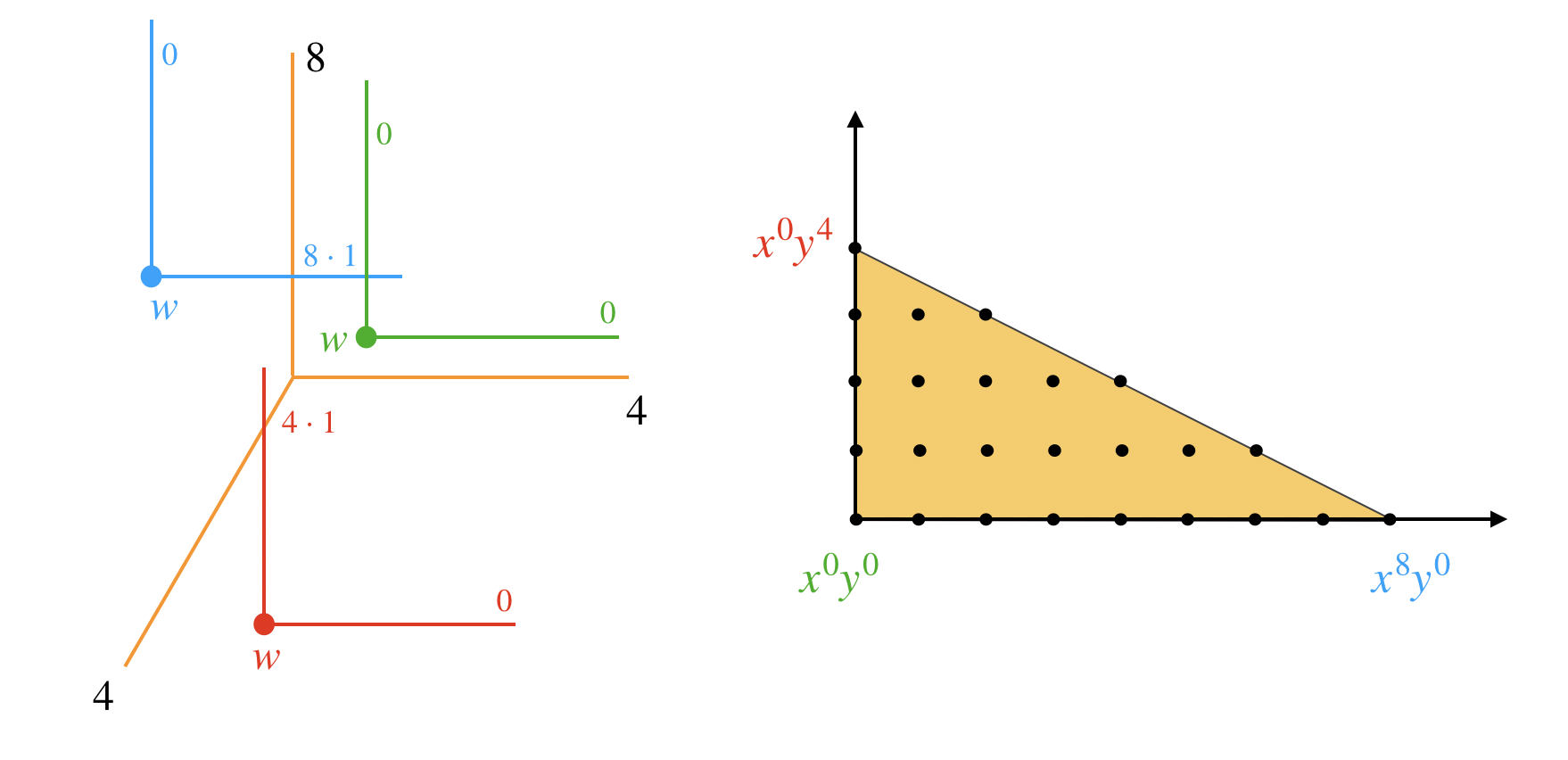}
\vspace{-0.15in}
\caption{Computing vertices of ${\cal N}(F)$ by intersecting ${\rm trop}(X)$ with $w + \mathbb{R}_{\geq 0} \cdot e_i$.}
\label{fig:raysintersect}
\end{figure}

The intersection multiplicity ${\rm IM}(w + \mathbb{R}_+ \cdot e_i, \sigma)$ is the \emph{lattice multiplicity} of the intersection of the ray $\mathbb{R} \cdot e_i$ with 
the hyperplane $\mathbb{L}_\sigma = \mathbb{R} \cdot \sigma$. This is the absolute value of the determinant of any $n \times n$ matrix whose columns are $e_i$ and a lattice basis for $\mathbb{L}_\sigma \cap \mathbb{Z}^n$. Algorithm \ref{alg: vertex_oracle} implements Theorem \ref{thm:troptonewt}.
It finds the vertex ${\cal N}(F)_w$ from the output of Algorithm \ref{alg:gentropimp} or \ref{alg:tropAdisc}.
Theorem \ref{thm:troptonewt} and Algorithm \ref{alg: vertex_oracle} are illustrated in Figure \ref{fig:raysintersect} for the  curve in Example \ref{ex:tropcurve}.

\begin{algorithm}
\caption{Compute vertex oracle from a tropical hypersurface}
\begin{algorithmic}[1]
\Procedure{getVertex(${\rm trop} (X \cap (\mathbb{C}^*)^n), \ w$)}{}
\State $v \gets 0$
\For{$i \in \{1, \dots, n\}$}
\For{ $(m_\sigma , \sigma) \in {\rm trop} (X) }$
	    \State $v \gets v +{\rm IM}(w+\RR_+ \cdot e_i, \sigma)\cdot m_\sigma \cdot e_i$
\EndFor
\EndFor
\State \Return \textit{$v$}
\EndProcedure
\end{algorithmic}
\label{alg: vertex_oracle}
\end{algorithm}

Algorithm \ref{alg: vertex_oracle} can be used to compute \emph{all} vertices of ${\cal N}(F)$. A naive approach applies the vertex oracle to many random vectors $w$. However, it is not clear how many 
$w$  would be needed, and which stop criterion to use. A deterministic way of constructing ${\cal N}(F)$ using a vertex oracle like Algorithm \ref{alg: vertex_oracle} was proposed by Huggins \cite{huggins}. Our implementation uses that.

\begin{example}
\label{ex: code: polytope from cycle}
The following {\tt Julia} code
computes a polytope from a tropical hypersurface: 
\begin{verbatim}
Delta = get_polytope_from_cycle(cone_list, weight_list)
\end{verbatim}
If the variables \texttt{cone\_list}, \texttt{weight\_list} are carried over from Example \ref{example: code: cycle frome parametrization}, then \texttt{Delta} is the yellow polytope
shown in Figure \ref{fig:raysintersect}.
For \texttt{cone\_list}, \texttt{weight\_list}  from Example \ref{example: code:  hyperdeterminant}, the polytope \texttt{Delta} is the Newton polytope of the hyperdeteminant $\Delta_A$
 from Example~\ref{ex: hyperdeterminant}.
\end{example}

\begin{remark} \label{rem:mixedfiber}
Under the assumptions of Section \ref{sec:2}, i.e., the $f_i$ are Laurent polynomials which are generic with respect to their Newton polytopes, ${\cal N}(F)$ is a \emph{mixed fiber polytope}. This was discovered independently by several authors \cite{emiris2007computing,esterov2008elimination,STY}. For instance, in Figure \ref{fig:polytopes}, ${\cal N}(F)$ is the mixed fiber polytope of $P_1$ and $P_2$. Our implementation of Huggins' algorithm  \cite{huggins} combined with Algorithm \ref{alg: vertex_oracle} provides a practical way of computing mixed fiber polytopes. This includes
the computation of fiber polytopes and secondary polytopes 
\cite[Section 3]{sturmfels2008tropical}.
\end{remark}

Once the Newton polytope ${\cal N}(F)$ of the defining equation $F= 0$ of $X$ is known, we can obtain its coefficients $c_a$ in \eqref{eq:ansatzF} using linear algebra. The set ${\cal B} = {\cal N}(F) \cap \mathbb{Z}^n$ is a superset of the monomial support ${\rm supp}(F)$ of $F$. It can be computed in \texttt{Oscar.jl} via the command \texttt{lattice\_points}. The interpolation method is most efficient when ${\cal B}$ is not much larger than ${\rm supp}(F)$, that is, few of the $c_a$ in \eqref{eq:ansatzF} are zero. 
The inclusion ${\cal B} \supseteq {\rm supp}(F)$ can be strict:
\begin{example}
Consider the map 
$f: \CC \rightarrow \CC^2$ given by
$f_1(t) =  a_1 t^4 + a_2 t\,$ and
$f_2(t) =   a_3 t^2 + a_4 t $ for generic complex numbers $a_1, a_2, a_3, a_4$.
Here the implicit polynomial $F(x,y)$ equals
$$  a_1^2   y^4 - 2 a_1 a_3^2   x y^2 + a_3^4   x^2 - 4 a_1 a_3 a_4^2   x y
+ 3 a_1 a_2 a_3 a_4   y^2 - a_4 (a_1 a_4^3-a_2 a_3^3)   x + a_2 (a_1 a_4^3-a_2 a_3^3)   y.
$$
Note that the term $y^3$ does not appear, in spite of it being
in $\mathcal{N}(F)$. This shows that
some lattice points in a predicted Newton polytope
may never appear with nonzero coefficient.
\end{example}
\begin{problem}
We propose to refine the observation in Remark \ref{rem:mixedfiber} by predicting the monomial support of $F$ from the monomial support of $f_1, \ldots, f_n$. That is, which lattice points in the mixed fiber polytope, other than its vertices, contribute to the implicit equation?
\end{problem}

For simplicity, we work with the superset ${\cal B} \supseteq {\rm supp}(F)$ and allow some coefficients to be zero. We identify a set ${\cal P}$ of $m$ points in $X$, so that the interpolation conditions $F(p) = 0$
for $ p \in {\cal P}$ uniquely determine $F$ (up to a constant factor). We obtain ${\cal P}\subset X$ by sending random points in $\mathbb{C}^d$ through the parametrization \eqref{eq:parametrization}. 
The unknown coefficients $c_a$ in \eqref{eq:ansatzF} form a vector $c = (c_a)_{a \in {\cal B}} \in \mathbb{C}^{\cal B}$.  
For each point $p \in \mathbb{C}^n$, let $p^{\cal B} = (p^a)_{a \in {\cal B}}$ be the vector of monomials corresponding to ${\cal B}$, evaluated at $p$. We interpret $p^{\cal B}$ as an element of the dual vector space $(\CC^{{\cal B}})^*$. With this set-up, $c$ is the unique vector (up to scaling) satisfying 
\[ p^{\cal B} \cdot c \,  = \, 0 \quad \text{ for all $p \in X$.}\]
If the sample points ${\cal P} \subset X$ are sufficiently random and $m \geq |{\cal B}|-1$, this is equivalent to 
\[ p^{\cal B} \cdot c \,  = \, 0, \quad \text{ for all $p \in {\cal P}$.}\]
The \emph{Vandermonde matrix} $M({\cal B},{\cal P})$ has the vectors $p^{\cal B}$ for  its rows,
where $ p \in {\cal P}$. It has size $m \times |{\cal B}|$ and, by the above discussion, the kernel of $M({\cal B}, {\cal P}) : \CC^{{\cal B}} \rightarrow \CC^m$ is spanned by $c$. 

Our problem is now reduced to the computation of the one-dimensional kernel of a Vandermonde matrix $M({\cal B}, {\cal P})$. Below, we will fix ${\cal B}$ and ${\cal P}$ and use the simpler notation $M = M({\cal B}, {\cal P})$, where there is no danger for confusion. When the parametrizing functions $f_i$ have coefficients in $\mathbb{Q}$, like in the case of $A$-discriminants in Section \ref{sec:4}, we can use ${\cal P} \subset \mathbb{Q}^n$. In particular, $M$ has rational entries, and its kernel can be computed in exact arithmetic.

\begin{example}
We now demonstrate our implementation of the above discussion
by computing the implicit equation $F$ from Example \ref{ex: tropical curve}.
The following code  computes the Vandermonde matrix $M({\cal B}, {\cal P})$ of size
 $24$ by $25$ with rational entries.
This is done by plugging $24$ random rational numbers into the parametrization \eqref{eq:parametrization}.
The functions \texttt{f1}, \texttt{f2}   are taken from Example \ref{example: code: cycle frome parametrization},
and the Newton polytope \texttt{Delta} $ = \mathcal{N}(F)$ was computed in Example \ref{ex: code: polytope from cycle}.
\begin{verbatim}
B = lattice_points(Delta)
n_samples = length(B)-1
P = sample([f1,f2], n_samples)
M_BP = get_vandermonde_matrix(B,P)
coeffs_F = nullspace(M_BP)[2]
\end{verbatim}
Up to scaling, \texttt{coeffs\_F} consists of the 25 coefficients of $F$. Some are shown in \eqref{eq:Fcurve}. 
\end{example}

 Often, in practical computations, the points $p \in {\cal P}$ are approximations of points on~$X$,
 so the entries of $M$ are finite precision floating point numbers. In that case, the task of computing $\ker M$ is one of \emph{numerical linear algebra}. This is not supported in the current version of \texttt{Oscar.jl}. The standard way to proceed using, for instance, the numerical linear algebra functionality in {\tt Julia}, is via the \emph{singular value decomposition} (SVD) of $M$. Alternatives include \emph{QR factorization with optimal pivoting} and \emph{iterative eigenvalue methods}.  
We refer to \cite[Section 5]{breiding2018learning} for
such numerical considerations and pointers to the relevant literature.

When $f$ is defined over $\QQ$, one might still want to use floating point computations for speed. Let $c_a$ be a nonzero entry of a generator $c$ for $\ker M$. The vector $c_a^{-1} c$ has rational entries. Its numerical approximation $\tilde{c}_a^{-1} \tilde{c}$ is contaminated by rounding errors. We approximate the entries of $\tilde{c}_a^{-1} \tilde{c}$ by rational numbers using the built in function \texttt{rationalize} in {\tt Julia}. This has an optional input \texttt{tol}, so that \texttt{rationalize(a,tol = e)} returns a rational number \texttt{q} which satisfies $|\mathtt{q} - \mathtt{a}| \leq \mathtt{e}$. A sensible choice for $\mathtt{tol}$ is $100 \cdot \tilde{c}_a^{-1} \cdot \varepsilon \cdot \sigma_1/\sigma_{|{\cal B}|-1}$. 

If symbolic computation is preferable to numerical methods, then one might solve
the linear equations over various finite fields and 
recover rational solutions via the Chinese remainder theorem.
This can be done in a computer algebra system.
Sometimes, one is only interested in a fixed finite field.
We illustrate the finite field computation in {\tt Oscar.jl}.

\begin{example}
We seek the $A$-discriminant for a matrix
whose entries are large integers:
$$ A \quad = \quad 
\begin{bmatrix}
1& 1& 1& 1& 1& 1\\
2& 3& 5& 7& 11& 13\\
 13& 8& 5& 3& 2& 1\\
\end{bmatrix}
$$
The following code finds that the 
 Newton polytope of $\Delta_A$ over $\QQ$
has dimension $3$ and f-vector $(12,18,8)$.
It terminated on a MacBook Pro with a 3,3 GHz Intel Core i5 processor within $120$ seconds.
The number of lattice points equals $2295$. In order to compute the coefficients of the $A$-discriminant, we must solve a linear system of $2294$ equations with
large integer coefficients.
We solve this over the field with $101$ elements instead:
\begin{verbatim}
A = [1 1 1 1 1 1; 2 3 5 7 11 13; 13 8 5 3 2 1];
cone_list, weight_list = get_trop_A_disc(A);
Delta = get_polytope_from_cycle(cone_list, weight_list);
@time mons, coeffs = compute_A_discriminant(A, Delta, GF(101));
\end{verbatim}
For  the same computation over the rational numbers, 
the machine ran out of memory. 
\end{example}

\smallskip

We close this section with a combinatorics problem that arises naturally from
Remark~\ref{rem:mixedfiber}.

\begin{problem} \label{prob:triangles}
Let $P_1,\ldots,P_n$ be polytopes in $\RR^{n-1}$ having
$v_1,\ldots,v_n$ vertices.
Give a sharp upper bound in terms of $v_1,\ldots,v_n$
for the number of vertices
of their mixed fiber polytope. In other words,
prove an Upper Bound Theorem for f-vectors
arising in tropical implicitization.
\end{problem}

\begin{example} \label{ex:triangles}
We illustrate Problem \ref{prob:triangles} for three triangles ($n=v_1=v_2=v_3=3$). After many runs for different random configurations, the following example is our current winner:
\begin{verbatim}
verts1 = [898 -614; -570 817; 892 -594]
verts2 = [-603 -481; -623 -127; -36 732]
verts3 = [-548 -864; -151 873; 800 -861]
(T1,T2,T3) = convex_hull.([verts1, verts2, verts3])
Delta = get_polytope_from_cycle(get_tropical_cycle([T1,T2,T3])...)
f_vec = f_vector(Delta)
\end{verbatim}
This code computes a mixed fiber polytope that has
 $25$ vertices, $49$ edges and $26$ facets.
 Can you find three triangles in $\RR^2$ whose mixed fiber polytope has more than $25$ vertices?
 \end{example}

\section{Higher codimension} \label{sec:5}

In this section we address the implicitization problem for
varieties $X$ that are not hypersurfaces. The role of the Newton polytope
$\mathcal{N}(F)$ of a polynomial $F$ will now be
played by the Chow polytope $\mathcal{C}(X)$.
We begin by reviewing some definitions from \cite{DS} and
 \cite[Chapter 6]{GKZ}.

Let $X$ be an irreducible projective variety of dimension $d$ in complex projective space $\PP^n$.
Suppose we are given  the tropical variety
${\rm trop}(X)$, a balanced fan of dimension $d$ in $\RR^{n+1}/\RR {\bf 1}$.
Our goal is to compute the {\em Chow form} ${\rm Ch}(X)$, which is a
hypersurface in the Grassmannian ${\rm Gr}(n-d-1,\PP^n)$. Its points
are the linear subspaces of dimension $n-d-1$ whose 
intersection with $X$ is non-empty.
We identify ${\rm Ch}(X)$ with its defining polynomial of degree ${\rm deg}(X)$ in
primal Pl\"ucker coordinates $p_{i_0 i_1 \cdots i_d}$, where
$1 \leq i_0 < i_1 < \cdots i_d \leq n$.  The $p_{i_0 i_1 \cdots i_d} $ are the maximal minors of
any $(d+1) \times (n+1)$
matrix whose kernel is the subspace.
The Chow form ${\rm Ch}(X)$ is only well-defined
up to the Pl\"ucker relations that vanish on ${\rm Gr}(n-d-1,\PP^n)$.
By \cite[Theorem 3.1.7]{AIT},  ${\rm Ch}(X)$ is a unique linear
combination of standard tableaux. In our computations, we always use that
standard representation for Chow forms.

The weight of the Pl\"ucker coordinate $p_{i_0 i_1 \cdots i_d}$
is the vector $e_{i_0} + e_{i_1} + \cdots + e_{i_d}$ in $\ZZ^n$,
and the weight of a Pl\"ucker monomial is the sum of the weights
of its variables, with multiplicity. By definition, the Chow polytope $\mathcal{C}(X)$
is the convex hull of the weights occurring in ${\rm Ch}(X)$.

Fink \cite{Fink} gave a combinatorial recipe for constructing the
weighted normal fan of the Chow polytope $\mathcal{C}(X)$ from the tropical variety ${\rm trop}(X)$.
Let $L_{n-d-1}$ denote the standard tropical linear space of dimension
$n-d-1$ in  $\RR^{n+1}/\RR {\bf 1}$. Its maximal cones are the
orthants spanned by $(n-d-1)$-tuples of unit vectors.
It is proved in \cite[Theorem 4.8]{Fink} that the weighted normal fan of $\mathcal{C}(X)$ is the
{\em stable sum} of ${\rm trop}(X)$ with the
negated linear space $- L_{n-d-1}$. The stable
sum is a dual operation to the stable intersection.
It always produces a balanced fan of expected dimension.
Hence ${\rm trop}(X) - L_{n-d-1}$ is a balanced
fan of codimension~$1$ in $\RR^{n+1}/\RR {\bf 1}$. 
Fink's result states that this is the
\emph{outer} normal fan of $\mathcal{C}(X)$. 

We can compute $\mathcal{C}(X)$ from ${\rm trop}(X) - L_{n-d-1}$ by the algorithm
for building Newton polytopes in Section~4, \emph{up to an integer translation}. Indeed, the normal fan of  $\mathcal{C}(X)$ and $\mathcal{C}(X) + \mathbf{t}$ are identical, for any $\mathbf{t} \in \mathbb{Z}^n$. Algorithm \ref{alg: vertex_oracle} finds vertices of ${\cal C}(X) + \mathbf{t}$, where $\mathbf{t}$ shifts ${\cal C}(X)$ so that it touches each coordinate hyperplane.  In previous examples, we had $\mathbf{t} = 0$. Indeed, if $F$ is irreducible,  then the polytope ${\cal N}(F)$ touches all coordinate hyperplanes. This is not true for the Chow polytope, as illustrated by the example below. Finding
the correct $\mathbf{t}$ is an interesting combinatorial problem which we plan to investigate in a future project.

\begin{example}[$d=1,n=3$] \label{ex:chow}
Let $X$ be the curve in $\CC^3$
which is given by the parametrization 
$$ x_1 = t(t-1)(t+1),\quad
x_2 =  t^2(t+1),\quad x_3 =  t^3(t-1). $$
The tropical curve is determined by
the orders of the coordinate functions at all
zeros and poles. Hence ${\rm trop}(X)$ is the
fan with four rays $(1,2,3),\, (1,1,0), \,(1,0,1)$ and $(-3,-3,-4)$.
We identify $X$ with its projective closure in $\PP^3$, obtained by
adding an extra coordinate $x_0$.
The tropical line $L_1$ is spanned by 
$e_0,e_1,e_2,e_3$, and we form the sum of 
${\rm trop}(X)$ with the negated line $-L_1$. This $2$-dimensional fan is the
normal fan of the Chow polytope $\mathcal{C}(X)$.

We implemented the stable sum using \texttt{Oscar.jl}, and obtain this fan as follows. 
\begin{verbatim}
cone_list = positive_hull.([[1, 1, 0], [1, 2, 3], [1,0,1], [-1, -1, -4//3]])
weight_list = ones(Int64, 4)
cone_list, weight_list = get_chow_fan(cone_list, weight_list)
\end{verbatim}
The output consists of 16 2-dimensional cones and their multiplicities. 
A translated version of the Chow polytope is obtained from this output as in the previous section: 
\begin{verbatim}
C_translated = get_polytope_from_cycle(cone_list, weight_list)
\end{verbatim}
This is a three-dimensional polytope touching all coordinate hyperplanes. It has vertices
$$ \begin{matrix} (0, 2, 3, 1), (0, 3, 1, 2), (0, 4, 1, 1), (1, 0, 4, 1), (1, 2, 3, 0), (1, 3, 0, 2), \\
     (1, 4, 0, 1), (1, 4, 1, 0),  (2, 0, 1, 3), (2, 0,4, 0), (2, 4, 0, 0), (3, 0, 0, 3). \end{matrix} $$
To identify the shift $\mathbf{t}$, we compare this to the Chow polytope. We obtain ${\cal C}(X)$ as the convex hull of the weights of the Pl\"ucker monomials in the Chow form ${\rm Ch}(X)$ of our curve:
  \begin{equation*} \begin{matrix}
p_{03}^4-p_{01}^3 p_{13} - 3 p_{01}^2 p_{02} p_{13} - 3 p_{01} p_{02}^2 p_{13}
-p_{02}^3 p_{13}+3 p_{01}^2 p_{03} p_{13} 
+9 p_{01} p_{02} p_{03} p_{13} + 6 p_{02}^2 p_{03} p_{13} \\
+ p_{01}
      p_{03}^2 p_{13}-5 p_{02} p_{03}^2 p_{13}+
2 p_{01}^2 p_{12} p_{13}+p_{01} p_{02} p_{12} p_{13} 
+2 p_{01}^2 p_{13}^2-2 p_{01} p_{02} p_{13}^2 +4 p_{02}^2p_{13}^2 \\
+p_{01} p_{03} p_{13}^2 - 4 p_{01} p_{12} p_{13}^2 
-p_{01}^3 p_{23}-3 p_{01}^2 p_{02} p_{23}-3 p_{01} p_{02}^2 p_{23} 
-p_{02}^3 p_{23}  + 4 p_{01}^2 p_{03} p_{23} \\ +11 p_{01} p_{02} p_{03} p_{23} 
+7 p_{02}^2 p_{03} p_{23} - 2 p_{01} p_{03}^2 p_{23} 
-10 p_{02} p_{03}^2 p_{23} 
+2 p_{03}^3 p_{23}+2 p_{01}^2 p_{12} p_{23}\\ +p_{01} p_{02} p_{12} p_{23}
+9 p_{01}^2 p_{13} p_{23}    -p_{01} p_{02} p_{13} p_{23} + 6 p_{02}^2 p_{13} p_{23}
      +2 p_{01} p_{03} p_{13} p_{23} 
      -2 p_{02} p_{03} p_{13} p_{23}\\ - 6 p_{01} p_{12} p_{13} p_{23}+2 p_{01} p_{13}^2 p_{23}+9 p_{01}^2 p_{23}^2
      +2 p_{01}  p_{02} p_{23}^2+2 p_{02}^2 p_{23}^2-4 p_{01} p_{03}p_{23}^2-2 p_{01} p_{12} p_{23}^2.
      \end{matrix}
\end{equation*}
We find that $\mathcal{C}(X)$ is the $3$-dimensional polytope with the following $12$ vertices:
$$ \begin{matrix} (1, 2, 3, 2), (1, 3, 1, 3), (1, 4, 1, 2), (2, 0, 4, 2), (2, 2, 3, 1), (2, 3, 0, 3), \\
     (2, 4, 0, 2), (2, 4, 1, 1),  (3, 0, 1, 4), (3, 0,4, 1), (3, 4, 0, 1), (4, 0, 0, 4). \end{matrix} $$
 We conclude that $\mathbf{t} = (-1, 0, 0, -1)$. For now, we apply this shift manually. 
     \end{example}

After listing all lattice points in $\mathcal{C}(X)$, 
we can compute the Chow form ${\rm Ch}(X)$ by interpolation.
This is done as follows.
For each lattice point $u$ in $\mathcal{C}(X)$ we 
list all standard Pl\"ucker monomials of weight $u$,
and form their linear combination with unknown coefficients.
Our ansatz is the sum of these 
$\ZZ^n$-homogeneous Pl\"ucker polynomials, with distinct unknown coefficients.
We generate random points on the Chow hypersurface as follows.
Pick a random point in $X$ and  a random linear space of dimension
$n-d-1$ through that point. We read off the Pl\"ucker coordinates
of that linear space and substitute them into the ansatz. 
Repeating this process many times gives the desired
 linear system of equations in the unknown coefficients. Up to scaling,
this system has a unique solution, namely the Chow form ${\rm Ch}(X)$.

\begin{example}
We use this strategy to recover the Chow form ${\rm Ch}(X)$ from Example \ref{ex:chow}. For each lattice point $u$ in the polytope ${\cal C}(X)$, we form the general linear 
combination of standard Pl\"ucker monomials of weight $u$. For instance, for $u= (2,2,2,2)$
this linear combination is
$$ \gamma_{u,1} \cdot p_{01}^2 p_{23}^2 \,+\, \gamma_{u,2} \cdot p_{01} p_{02} p_{13} p_{23}
\,+\, \gamma_{u,3} \cdot p_{02}^2 p_{13}^2. $$
Our ansatz for  the Chow form ${\rm Ch}(X)$ is the sum of these expressions over all $u \in \mathcal{C}(X) \cap \ZZ^4$.

We sample from the Chow hypersurface by picking random matrices of the form
$$ \begin{bmatrix}\,\alpha_0  & \alpha_1 & \alpha_2  & \alpha_3 \\ \,
1 &  t(t-1)(t+1) &  t^2(t+1) & t^3(t-1) \end{bmatrix}. $$
The $2 \times 2$ minors of this matrix are the dual Pl\"ucker coordinates of the 
corresponding sample point
in ${\rm Ch}(X) \subset {\rm Gr}(1,\PP^3)$.
We read off  its primal Pl\"ucker coordinates as follows:
$$ \small  \begin{matrix}
p_{01} =  (t^4-t^3)\alpha_2+(t^3+t^2)\alpha_3,  &
p_{02} =  (t^3-t^4)\alpha_1+(t^3-t)\alpha_3,  & 
p_{03} =  (t^3+t^2)\alpha_1+(t-t^3)\alpha_2,  \\
p_{12} = (t^4-t^3)\alpha_0-\alpha_3,\, & 
p_{13} = (t^3+t^2)\alpha_0+\alpha_2,\, & 
p_{23} = (t^3-t)\alpha_0-\alpha_1. 
\end{matrix}
$$
We substitute many such sample points into the ansatz, and we solve
the resulting system of linear equations for the unknown coefficients $\gamma_{u,i}$.
The output is the desired Chow form.
This yields defining equations for $X$ by setting
$p_{ij} = \alpha_i x_j - \alpha_j x_i$ for any  $\alpha_0,\ldots,\alpha_3 \in \QQ$.
\end{example}

In this article, we discussed an application of 
tropical geometry to computer algebra, namely 
implicitization with tropical preprocessing. It was shown
that {\tt Oscar.jl} provides excellent capabilities for
performing tropical implicitization in practice. Our implementation
in {\tt Oscar.jl} realizes the vision in  \cite{sturmfels2008tropical} and
 fulfils the promise made by  {\tt TrIm}.
In Section~5 we ventured into
a setting where the desired hypersurface is not in an affine or projective space,
but inside a Grassmannian.
This suggests yet one more
problem for future research.

\begin{problem}
Many applications lead to interesting subvarieties of Grassmannians.
For instance, in computer vision, certain cameras are represented 
by curves and surfaces in ${\rm Gr}(1,\PP^3)$. 
Their tropicalizations lie inside the tropical Grassmannian,
and their cohomology classes are computed  by Schubert calculus.
It would be desirable to develop tropical implicitization 
in the setting when the ambient spaces are Grassmannians, or even
 flag varieties.
\end{problem}

\bigskip

\bibliographystyle{abbrv}

\begin{thebibliography}{10}

\bibitem{breiding2018learning}
P.~Breiding, S.~Kali{\v{s}}nik, B.~Sturmfels, and M.~Weinstein:
\newblock Learning algebraic varieties from samples.
\newblock {\em Revista Matem{\'a}tica Complutense}, 31:545--593, 2018.

\bibitem{sandra}
E.~Cattani, M.~Cueto, A.~Dickenstein, S.~Di Rocco, and B.~Sturmfels:
\newblock Mixed discriminants. 
\newblock {\em Mathematische Zeitschrift}, 274: 761--778, 2013.

\bibitem{corless2001numerical}
R.~M. Corless, M.~W. Giesbrecht, I.~S. Kotsireas, and S.~M. Watt.
\newblock Numerical implicitization of parametric hypersurfaces with linear
  algebra:
\newblock In {\em Artificial Intelligence and Symbolic Computation:
  International Conference AISC 2000 Madrid, Spain, July 17--19, 2000 Revised
  Papers 5}, pages 174--183. Springer, 2001.

\bibitem{cox1997ideals}
D.~Cox, J.~Little, and D.~O’Shea:
\newblock {\em Ideals, Varieties, and Algorithms}.
\newblock Undergraduate Texts in Mathematics. Springer-Verlag, New York, 1997.

\bibitem{DS}
J.~Dalbec and B.~Sturmfels: Introduction to Chow forms,
{\em Invariant Methods in Discrete and Computational Geometry} (Cura\c{c}ao, 1994), 37--58, Kluwer Acad. Publ., Dordrecht, 1995.

\bibitem{dickenstein2007tropical}
A.~Dickenstein, E.~Feichtner, and B.~Sturmfels:
\newblock Tropical discriminants.
\newblock {\em Journal of the American Mathematical Society}, 20(4):1111--1133,
  2007.

 \bibitem{emiris2013implicitization}
 I.~Z. Emiris, T.~Kalinka, C.~Konaxis, and T.~L. Ba:
 \newblock Implicitization of curves and (hyper) surfaces using predicted  support.
 \newblock {\em Theoretical Computer Science}, 479:81--98, 2013.

\bibitem{emiris2007computing}
I.~Z. Emiris, C.~Konaxis, and L.~Palios:
\newblock Computing the {N}ewton polytope of specialized resultants.
\newblock Presented at MEGA 2007 (Effective methods in algebraic geometry).

\bibitem{esterov2008elimination}
A.~Esterov and A.~Khovanskii:
\newblock Elimination theory and Newton polytopes.
\newblock {\em Functional Analysis and Other Mathematics}, 2(1):45--71, 2008.

\bibitem{Fink}
A.~Fink: 
\newblock Tropical cycles and Chow polytopes.
\newblock {\em Beitr.~Algebra Geom.}~54:13-40, 2013.

\bibitem{GKZ}
I.~Gel'fand, M.~Kapranov, and A.~Zelevinsky: {\em Discriminants, Resultants, and Multidimensional Determinants}, 
 Birkh\"auser, Boston, 1994.

\bibitem{huggins}
P.~Huggins:
\newblock i{B}4e: A software framework for parametrizing specialized {LP} problems.
\newblock In A.~Iglesias and N.~Takayama, editors, {\em Mathematical Software -
  ICMS 2006}, pages 245--247, Springer Verlag, Berlin-Heidelberg, 2006. 
  

\bibitem{maclagan2021introduction}
D.~Maclagan and B.~Sturmfels:
\newblock {\em Introduction to Tropical Geometry}, Graduate Texts in Mathematics, volume 161.
\newblock American Mathematical Society, 2021.

\bibitem{AIT}
B.~Sturmfels: {\em Algorithms in Invariant Theory}, Texts and Monographs in Symbolic Computation, Springer-Verlag, Vienna, 1993.

\bibitem{ST}
B.~Sturmfels and J.~Tevelev:
\newblock Elimination theory for tropical varieties.
\newblock {\em Mathematical Research Letters}, 15:543--562, 2008.

\bibitem{STY}
B.~Sturmfels, J.~Tevelev, and J.~Yu:
\newblock The {N}ewton polytope of the implicit equation.
\newblock {\em Moscow Mathematical Journal}, 7:327--346, 2007.

\bibitem{sturmfels2008tropical}
B.~Sturmfels and J.~Yu:
\newblock Tropical implicitization and mixed fiber polytopes.
\newblock {\em Software for Algebraic Geometry}, 111--131, 
IMA Vol. Math. Appl., 148, Springer, New York, 2008. 



\end{thebibliography}

		\bigskip
		\bigskip
		
		\noindent {\bf Authors' addresses:}
		
		\smallskip
			
		\noindent Kemal Rose, MPI-MiS Leipzig
		\hfill \url{kemal.rose@mis.mpg.de}
		
		\noindent  Bernd Sturmfels, MPI-MiS Leipzig \hfill \url{bernd@mis.mpg.de}

		\noindent  Simon Telen, MPI-MiS Leipzig \hfill \url{simon.telen@mis.mpg.de}

\end{document}